\documentclass[a4paper,12pt,twoside]{article}
\usepackage[italian]{babel}
\pdfoutput=1
\usepackage{amsmath}
\usepackage{amssymb}
\usepackage{latexsym}
\usepackage[mathscr]{euscript}
\usepackage{verbatim}
\usepackage{graphicx}
\DeclareGraphicsExtensions{.png}
\begin{document}
\vspace*{10mm}
\thispagestyle{empty} \Large \noindent
\textbf{Il Fattore di \vspace{20mm} Sylvester} \\
\large
\textbf{Donato Saeli $\;\cdot\;$ Maurizio \vspace{43mm} Spano} \\
\small
\pmb{\textsc{Abstract }} Sylvester factor, an essential part of the asymptotic formula of Hardy and Littlewood which is the \textit{extended} Goldbach conjecture, regarded as strongly multiplicative arithmetic function, has several remarkable \vspace{2mm} properties.

Il fattore di Sylvester, parte essenziale della formula asintotica di Hardy e Littlewood che costituisce la congettura\, \textit{estesa}\, di Goldbach, riguardato come funzione aritmetica fortemente moltiplicativa, presenta diverse propriet\`a \vspace{3mm} significative. \\
\pmb{\textsc{Keywords }} Goldbach's extended conjecture, Sylvester factor. $\;\cdot\;$ Strongly multiplicative function. $\;\cdot\;$ Convolution \vspace{3mm} inverse. \\
\pmb{\textsc{MSC }} \vspace{\stretch{1}} 11P32 \\
\rule[.3mm]{90mm}{.3mm} \vspace{1.5mm} \\
Donato Saeli \\
Via Giovanni XXIII, 29 - 85100 Potenza (PZ), Italy \\
Tel.: +39-0971-51280 \\
Email: \vspace{2.5mm} donato.saeli@gmail.com \\
Maurizio Spano \\
Via Rocco Scotellaro, 19 - 75019 Tricarico (MT), Italy \\
Tel.: +39-348 998 5205

\setcounter{footnote}{1}
\noindent \normalsize
\textbf{1\, \vspace{2mm} Introduzione}

Consideriamo la funzione aritmetica $\, g(n)\, , \,$ che associa ad $\, n \,$ il numero delle coppie ordinate $\, (p,q) \,$ di numeri p\underline{rimi} \underline{dis}p\underline{ari} tali che $\, {p + q = 2n}; \,$ il grafico di $\, g(n) \,$ in $\, [3,N], \,$ con $\, N \,$  sufficientemente grande, ap-\\
pare come \vspace{2mm} una cometa.\,\footnote{ \ La \textit{cometa di Goldbach}: l'affermazione \ $g(n) > 0 \ \ per \ ogni \ \ n > 2 \ $ equivale \\
alla congettura di Goldbach: \textit{Ogni numero pari non inferiore a quattro \`e somma \\
di due primi}.}

La congettura\, \textit{estesa}\, di Goldbach, formulata nel 1922 da Hardy \\
e Littlewood afferma che
$$g(n) \sim h(n) := \dfrac{4 \, c \, n}{(\lg n)^2} \prod_{
	\begin{subarray}{c}
	p \, | \, n \\
	p \neq 2
	\end{subarray}}
\dfrac{p-1}{p-2} \; , \hspace{21mm} (1)$$
dove Il prodotto $\, \displaystyle \prod_{
	\begin{subarray}{c}
	p \, | \, n \\
	p \neq 2
	\end{subarray}}
\dfrac{p-1}{p-2} \,$ s'intende esteso a tutti i numeri primi \vspace{1mm} dispari che dividono $\, n$ e si pone uguale a 1 se \`e privo di fattori, cio\`e se \vspace{3mm} $\, n = 2^k \ (k \geq 0).$ \\
La costante $\, c \,$ \`e il valore del prodotto (infinito), esteso a tutti i primi dispari: \\
$$ \prod_{p \neq 2} \dfrac{p(p-2)}{(p-1)^2} = \lim_{n\rightarrow\infty} \bigg[\prod_{2 < p < n} \Big( 1 - \dfrac{1}{(p-1)^2} \Big) \bigg] = \vspace{2mm} 0,6601618. \dots$$

Sylvester {[}S{]} nel 1871, aveva proposto una formula equivalente \vspace{1mm} alla \\
$\, g(n) \sim 2 \, e^{-\gamma} \, h(n); $\,\footnote{ \ $2 \, e^{-\gamma} = 1,1229. \dots$}  \ ma la  \vspace{2mm} (1) \\
``...  \`e la sola formula di questa sorta che pu\`o essere corretta, cosicch\'e la formula di Sylvester \`e errata. \ Ma Sylvester \`e stato il primo ad identificare il fattore
$$\prod_{
	\begin{subarray}{c}
	p \, | \, n \\
	p \neq 2
	\end{subarray}} \vspace{1mm}
\dfrac{p-1}{p-2}$$
a cui sono dovute le \textit{irregolarit\`a} della $\, h(n). \,$ Non vi sono indicazioni sufficienti per mostrare come sia stato condotto al suo risultato. ...'' {[}HL{]} (pp.\! 32, \vspace{2mm} 33).

La (1) pu\`o essere riguardata in vari modi {[}Smd{]} (pp.\! 2-4), molto suggestivo \`e il seguente:
$$\mathscr{S}(n) :=
\prod_{
	\begin{subarray}{c}
	p \, | \, n \\
	p \neq 2
	\end{subarray}} \vspace{1mm}
\dfrac{p-1}{p-2} \sim G(n) := \dfrac{(\lg n)^2}{4 \, c \, n}\, g(n). \hspace{18mm} (2)$$
Se si effettua il confronto fra i due membri della (2), il risultato inizialmente irrilevante (fig. 1), diventa piuttosto interessante al crescere di n (fig. 2).
\begin{center}
	\includegraphics[width=1\textwidth]{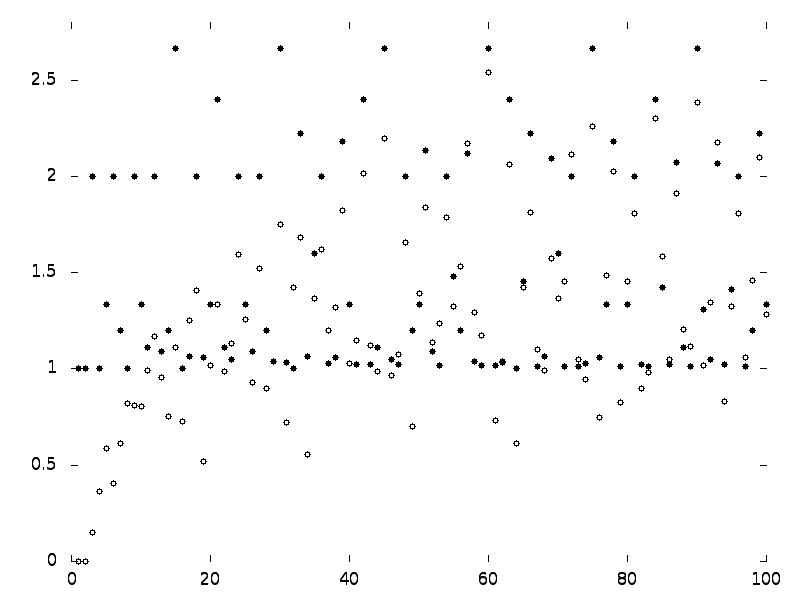} \\
	Fig. 1 \quad \includegraphics[scale =.5]{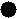} \ $\mathscr{S}(n),$ \ \ \includegraphics[scale =.5]{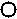} \ \vspace{6mm} $\, G(n)$ \\
	\includegraphics[width=1\textwidth]{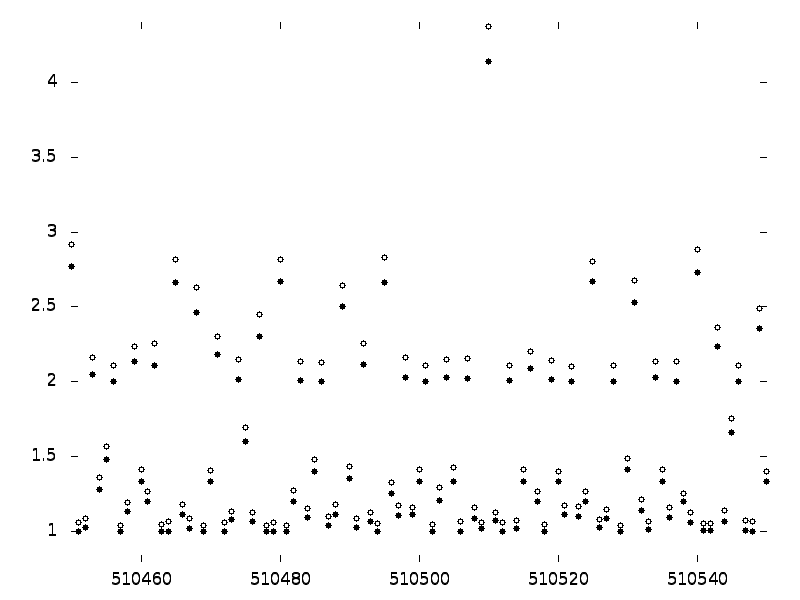} \\
	Fig. 2 \quad \includegraphics[scale =.5]{tnp.png} \ $\mathscr{S}(n),$ \ \ \includegraphics[scale =.5]{tnv.png} \ $\, G(n)$ \\
\end{center}
Gli autori hanno verificato che \`e $\, \mathscr{S}(n) < G(n) \, , \,$ per $\, 72.065 \leq n \leq 2.000.000\, . \ $ Sarebbe interessante stabilire se e per quale valore di $\, n \ (> 2 \cdot 10^6) \,$ risulti $\, \mathscr{S}(n) > \vspace{2mm} G(n) \, .$

Vogliamo richiamare infine l'attenzione su un'altra relazione nella quale il \textit{fattore di Sylvester} riveste un ruolo essenziale. Se $\, m \in \mathbb{N}, \,$ indichiamo con $\, \mathbb{Z}_m \,$ l'anello delle classi di resto modulo $\, m \,$ e con $\, \mathbb{Z}_m^{^*} \,$ il gruppo delle classi prime con $\, m. \,$ \`E noto che
$$ s_m^*(n) = \sharp\, \{(\bar{r},\bar{s}) \in \mathbb{Z}_m^{^*} \times \mathbb{Z}_m^{^*}\; |\; \bar{r}+\bar{s} = \bar{n}\} = m \prod_{
	\begin{subarray}{n}
	p \, | \, m \\
	p \, | \, n
	\end{subarray}
}\bigg(1-\dfrac{1}{p} \bigg) \prod_{
	\begin{subarray}{n}
	p \, | \, m \\
	p \, \nmid \, n
	\end{subarray}
} \bigg(1-\dfrac{2}{p} \bigg) $$
{[}Dm{]}, {[}Sj{]}; \ se scegliamo $\, m = 2 \cdot q_1 \cdots q_t, \,$ con $\, q_1 < q_2 < \cdots < q_t \,$ \vspace{2mm} numeri primi dispari e per $\, n \in \mathbb{N} \,$  poniamo
$\displaystyle \alpha_i = \left\{ \begin{array}{l}
1 \ \; \textrm{se} \ \; q_i | n \\
0 \ \; \textrm{se} \ \; q_i \nmid n  \\
\end{array} \right. \ \ \textrm{per} \ \ i = 1, \, \dots, \, \vspace{2mm} t,$
allora $\, d = (2n, \, m) = 2 \cdot q_1^{\alpha_1} \cdots p_t^{\alpha_t} \,$ e quindi anche
$$ s_m^*(2n) = m \prod_{p \, | \, d}\bigg(1-\dfrac{1}{p} \bigg) \prod_{p \, | \, \frac{m}{d}} \bigg(1-\dfrac{2}{p} \bigg) = \prod_{i=1}^t (q_i-1 )^{\alpha_i} (q_i-2)^{1-\alpha_i} $$
$$ = \prod_{i=1}^t \bigg( \dfrac{q_i-1}{q_i-2} \bigg)^{\alpha_i} \prod_{i=1}^t (q_i-2). \ \qquad \textrm{In altri termini:}$$ \vspace{.5mm}
$$s_m^*(2n) = \vspace*{6mm} \mathscr{S}(d) \cdot s_m^*(2).$$
\textbf{2\, Osservazioni sulle funzioni fortemente \vspace{2mm} moltiplicative}

Una funzione aritmetica $\, f(n) \,$ si dice \textit{fortemente moltiplicativa} se \`e moltiplicativa e quali che siano $\, p \,$ primo, $\, k \,$ intero positivo, \`e \vspace{3mm} $\, f(p^k) = f(p).$ \\
Esempi di funzioni fortemente moltiplicative \vspace{2mm} sono: \\
la funzione $ \quad \displaystyle \overline{\varphi}(n) = \prod_{p \, | \, n} \dfrac{p-1}{p} = \dfrac{\varphi(n)}{n} $ \qquad ($ \varphi(n) \,$ indicatore di \vspace{-1mm} Eulero) \\
e giust'appunto il fattore di Silvester $ \quad \displaystyle \mathscr{S}(n) = \prod_{
	\begin{subarray}{c}
	p \, | \, n \\
	p \neq 2
	\end{subarray}} \vspace{3mm}
\dfrac{p-1}{p-2} \qquad $ (fig. 3 e 4). \\
\textsc{Proposizione 2.1} - \it Se $\, f \,$ e $\, g \,$ sono funzioni aritmetiche fortemente moltiplicative, $\, p \,$ un primo e $\, k \,$ un intero positivo, allora \\
$\, (f \ast g)(p^k) = \displaystyle \sum_{i=0}^{k} f(p^i) g(p^{k-i}) = f(p) + g(p) + (k-1) f(p) g(p) $ \\
e per $\, k \geq 2, \qquad f^{-1}(p^k) = (-1)^k f(p) \big[ f(p) -1 \big]^{k-1}\, . $ \vspace{3mm} \rm \hspace{\stretch{1}} $\bullet$ \\
\begin{center}
	\includegraphics[width=1\textwidth]{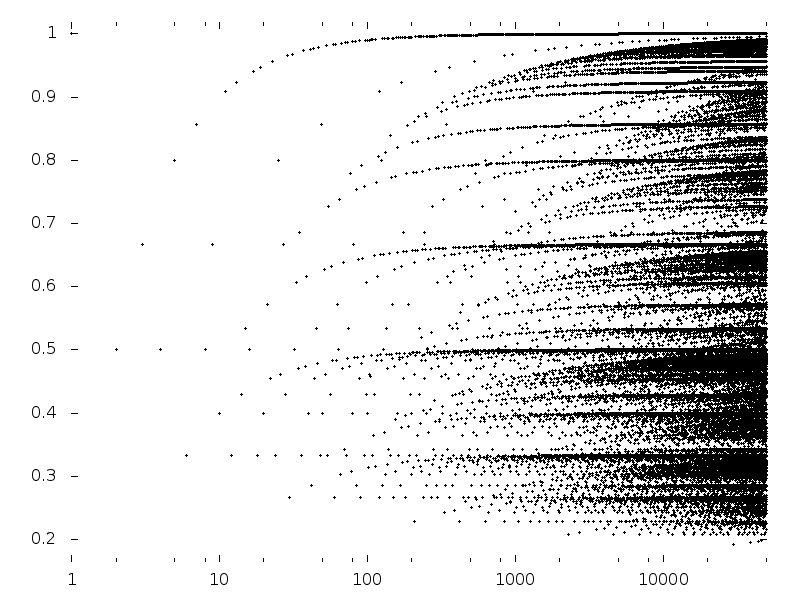} \\
	Fig. 3 \qquad \vspace{6mm} $ \overline{\varphi}(n) = \dfrac{\varphi(n)}{n} $ \\
	\includegraphics[width=1\textwidth]{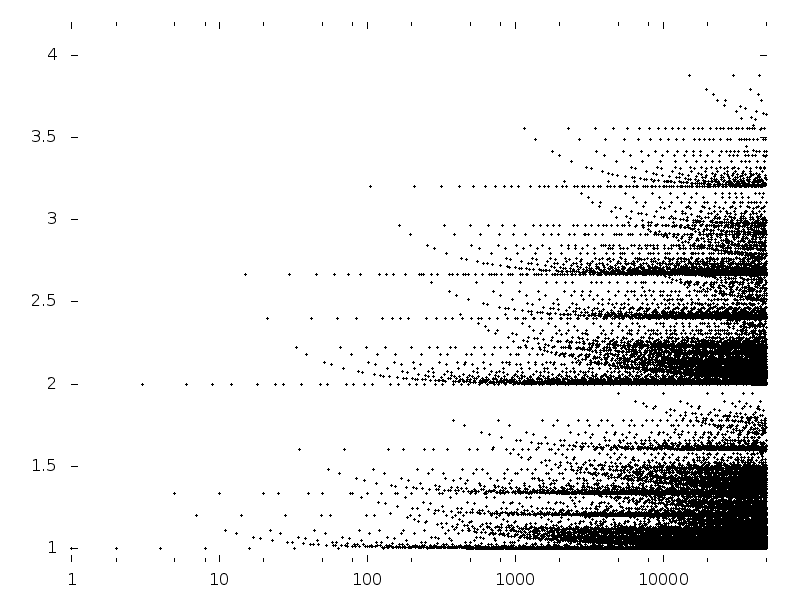} \\
	Fig. 4 \qquad $\mathscr{S}(n)$ \\
\end{center}
\textsc{Dimostrazione.} \ (Per induzione su \vspace{2mm} $\, k$). \\
$ f^{-1}(p^2) = \vspace{1mm} -f(p) f^{-1}(p) - f(p^2) f^{-1}(1) =  f(p) \big[ f(p) -1 \big]. $ \\
$ f^{-1}(p^{k+1}) = \vspace{1mm} \displaystyle - \sum_{i=1}^{k+1} f(p^i) f^{-1}(p^{k-i+1}) $ \\
$ \displaystyle = -f(p) \Big\{ -f(p) + 1 + \sum_{i=1}^{k-1} f^{-1}(p^{k-i+1}) \Big\} $ \\
$ \displaystyle = -f(p) \Big\{ -f(p) + 1 + \sum_{i=1}^{k-1} (-1)^{k-i+1} f(p) \big[ f(p) -1 \big]^{k-i} \Big\} $ \\
$ \displaystyle = f(p) \big[ f(p) -1 \big] \Big\{ 1 + \sum_{i=1}^{k-1} (-1)^{k-i} f(p) \big[ f(p) -1 \big]^{k-i-1} \Big\} $ \\
$ \displaystyle = f(p) \big[ f(p) -1 \big] \sum_{i=1}^k f^{-1}(p^{k-i}) = \vspace{1mm} \big[ f(p) -1 \big] \sum_{i=1}^k f(p^i) f^{-1}(p^{k-i}) $ \\
$ = - \big[ f(p) -1 \big] f^{-1}(p^k) = \vspace{2mm} - \big[ f(p) -1 \big] (-1)^k f(p) \big[ f(p) -1 \big]^{k-1}$ \\
$ = (-1)^{k+1} f(p) \big[ f(p) -1 \big]^k$. \vspace{4mm} \hspace{\stretch{1}} $\vartriangle$ \\
Cos\`{\i} \vspace{2mm} abbiamo: \\
$\, ( \overline{\varphi} \ast \mathscr{S} )(p^k) = (k+1) \Big( 1 + \dfrac{1}{p(p-2)} \Big)\, ,$ \ se $\, p \neq 2 \, ; \qquad ( \overline{\varphi} \ast \mathscr{S} )(2^k) = \vspace{1mm} 1 + \dfrac{k}{2} \, .$ \\
$\, \bar{\varphi}^{-1}(p^k) = \vspace{2mm} \dfrac{1-p}{p^k}\, .$ \\
$\, \mathscr{S}^{-1}(p^k) = \dfrac{1-p}{(p-2)^k}\, ,$ \ se $\, p \neq 2 \, ; \qquad \mathscr{S}^{-1}(2) = -1 \, , \ \ \mathscr{S}^{-1}(2^k) = 0 \, ,$ se \vspace{4mm} $\, k \geq 2 \, .$ \\
\textsc{Proposizione 2.2} - \it Sia $\, f \,$ una funzione aritmetica fortemente moltiplicativa. Per ogni intero $\, m > 1 \, , \,$  posto $\, y = f(m) \,$ e $\, f^-(y) = \{ n \in \mathbb{N}\, |\, f(n) = y \}\, , \,$ si ha$:$ \vspace{2mm} $\, | f^-(y) | = \aleph_0 \, .$ \rm \hspace{\stretch{1}} $\bullet$ \\
Notiamo che $\, \bar{\varphi}^-(1) = \{1\}\, . \,$ mentre \vspace{3mm} $\, | \mathscr{S}^-(1) | = \aleph_0 \, .$ \\
\textbf{3\, Propriet\`a  particolari delle funzioni} $\, \boldsymbol{\bar{\varphi}} \,$ \textbf{e} \vspace{2mm} $\, \boldsymbol{\mathscr{S}}$

Se indichiamo con $\, \{p_n\} \,$ la successione dei numeri primi in ordine \vspace{2mm} crescente, \`e chiaro che la successione $\, \{\overline{\varphi}(p_n)\} \,$ \`e monotona crescente, \ la $\, \{\mathscr{S}(p_n)\} \,$ \vspace{2mm} \`e, \\
per $\, n > 1, \,$ monotona decrescente \vspace{2mm} e \ $\, \displaystyle \lim_{n\rightarrow\infty} \overline{\varphi}(p_n) = \lim_{n\rightarrow\infty} \mathscr{S}(p_n) = 1.$ \\
\textsc{Proposizione 3.1} - \it $\, \mathscr{R}(\mathscr{S}) \,$ e $\, \mathscr{R}(\overline{\varphi}) \,$ sono entrambi insiemi perfetti e totalmente \vspace{1mm} sconnessi. \rm \hspace{\stretch{1}} $\bullet$ \\
\textsc{Dimostrazione.} \ Infatti $\, \forall\, n \in \mathbb{N}, \ \mathscr{S}(n) \in \mathbb{Q}\, ; \quad $ e se $\, p > n \,$ \`e un primo dispari, allora $\, \mathscr{S}(np) \neq \mathscr{S}(n) \,$ e $\, |\mathscr{S}(np) - \mathscr{S}(n)|=|\mathscr{S}(n)||\mathscr{S}(p) - 1| < \varepsilon\, , \,$ \\
con $\, \varepsilon > 0 \,$ arbitrario, posto sia anche $\, p > M_\varepsilon \, .$ \ Considerazioni simili valgono per la \vspace{1mm} $\, \overline{\varphi}\, .$ \hspace{\stretch{1}} $\vartriangle$ \\
Definiamo la successione $\, \{P_n\} \,$ come \vspace{2mm} segue: \\
\hspace*{9mm} \vspace{2mm} $ P_1 = p_1, \ \ P_2 = p_1 \cdot p_2, \ \ \dots, \ \ P_n = p_1 \cdot p_2 \cdots p_n, \ \ \cdots \ .$ \\
Abbiamo: $\, \displaystyle \overline{\varphi}(P_n) = \prod_{i=1}^{n}\bigg(1-\dfrac{1}{p_i}\bigg) \quad $ e \vspace{2mm} quindi \\
$\, \displaystyle \lim_{n\rightarrow\infty} \overline{\varphi}(P_n) = \prod_{i=1}^{\infty} \bigg( 1 - \dfrac{1}{p_i} \bigg) = 0\, ; \ \ $ dato che \vspace{2mm} $\ \displaystyle \sum_{i=1}^{\infty} \dfrac{1}{p_i} = +\infty\, . $ \\
Per lo stesso motivo \vspace{3mm} $\ \displaystyle \lim_{n\rightarrow\infty} \mathscr{S}(P_n) = \prod_{i=2}^{\infty} \bigg( 1 + \dfrac{1}{p_i-2} \bigg) = +\infty\, . $ \\
\textsc{Proposizione 3.2} - \it Se $\, m < P_n, \,$ allora $\, \overline{\varphi}(m) > \overline{\varphi}(P_n)\, .$ \\
Se $\, m < P_n \,$ e $\, 2m \neq P_n, \,$ allora $\, \mathscr{S}(m) < \mathscr{S}(P_n)\, .$ \rm \vspace{2mm} \hspace{\stretch{1}} $\bullet$ \\
\textsc{Dimostrazione.} \ Infatti se $\, m < P_n \,$ e $\, q_1^{\alpha_1} \cdots q_\nu^{\alpha_\nu} = m \,$ ne \`e la decom- \\
posizione canonica in fattori primi, allora deve essere $\, \nu < n \,$ e $\, q_i \geq p_i \, ,$ \\
per \vspace{.5mm} $\, i = 1, \dots, \nu \, .$ \\
Ne segue \vspace{2mm} $\ \displaystyle \overline{\varphi}(m) = \prod_{i=1}^{\nu} \bigg( 1 - \dfrac{1}{q_i} \bigg) \geq \prod_{i=1}^{\nu} \bigg( 1 - \dfrac{1}{p_i} \bigg) > \prod_{i=1}^{\nu+1} \bigg( 1 - \dfrac{1}{p_i} \bigg) \geq \overline{\varphi}(P_n)\, .$ \\
Se $\, m \,$ \`e pari, \vspace{2mm} abbiamo \\
$\ \displaystyle \mathscr{S}(m) = \prod_{i=2}^{\nu} \bigg( 1 + \dfrac{1}{q_i-2} \bigg) \leq \prod_{i=2}^{\nu} \bigg( 1 + \dfrac{1}{p_i-2} \bigg) < \prod_{i=2}^{\nu+1} \bigg( 1 + \dfrac{1}{p_i-2} \bigg) \leq \vspace{2.5mm} \mathscr{S}(P_n)\, .$ \\
Se $\, m \,$ \`e dispari e $\, \nu < n - 1\, , \,$ \vspace{2mm} abbiamo $\, q_i \geq p_{i+1} \, , \,$ per $\, i = 1, \dots, \nu \,$  e \\
$\ \displaystyle \mathscr{S}(m) = \prod_{i=1}^{\nu} \bigg( 1 + \dfrac{1}{q_i-2} \bigg) \leq \prod_{i=2}^{\nu+1} \bigg( 1 + \dfrac{1}{p_i-2} \bigg) < \prod_{i=2}^{\nu+2} \bigg( 1 + \dfrac{1}{p_i-2} \bigg) \leq \vspace{2.5mm} \mathscr{S}(P_n)\, .$ \\
Se $\, m \,$ \`e dispari e $\, \nu = n - 1\, , \,$ \vspace{1mm} allora deve essere necessariamente $\, q_i = p_{i+1} \, , \,$ per $\, i = 1, \dots, n - 1 \,$  e \vspace{5mm} $\ \mathscr{S}(m) = \mathscr{S}(P_n)\, .$ \hspace{\stretch{1}} $\vartriangle$ \\
\large
\textbf{Bibliografia}\vspace{3mm} \\
\normalsize
\begin{tabular}{ll}
{[}Dm{]} & \textsc{M. Deaconescu}, \textit{Adding units mod n}, \\
& Elem. Math., 55(2000), 123-127. \\
{[}HL{]} & \textsc{G. H. Hardy - J. E. Littlewood}, \ \textit{Some problems} \\
& \textit{of `partitio numerorum' III: on the expression of a number} \\
& \textit{as a sum of primes}, Acta Math., 44(1922), 1-70. \\
{[}S{]} & \textsc{J. J. Sylvester}, \textit{On the partition of an even number into} \\
& \textit{two primes}, Proc. London Math. Soc. Ser.\,I, 4(1871), 4-6. \\
{[}Sj{]} & \textsc{J. W. Sander}, \textit{On the addition of units and nonunits mod m}, \\
& Journal of Number Theory, 129(2009), 2260-2266. \\
{[}Smd{]} & \textsc{D. Saeli - M. Spano}, \textit{La cometa di Goldbach e ... le altre}, \\
& Lecture Notes of \textit{Seminario interdisciplinare di Matematica}, \\
& Vol. 10(2011), 45-57.
\end{tabular}

\end{document}